\documentclass[12pt,thmsa]{article}%
\usepackage{sw20bams}%
\usepackage{amsmath}%
\usepackage{amsfonts}%
\usepackage{amssymb}%
\usepackage{graphicx}
\providecommand{\U}[1]{\protect\rule{.1in}{.1in}}
\begin{document}

\author{Steven Finch}
\title{Two Sample Covariances from a Trivariate Normal Distribution}
\date{December 17, 2015}
\maketitle

\begin{abstract}
The joint distribution of two off-diagonal Wishart matrix elements was useful
in recent work on geometric probability \cite{F1}. Not finding such formulas
in the literature, we report these here.

\end{abstract}

\footnotetext{Copyright \copyright \ 2010, 2015 by Steven R. Finch. All rights
reserved.}Let $(A,B,C)$ be trivariate normally distributed with known mean
$(0,0,0)$ and covariance matrix
\[%
\begin{array}
[c]{ccccc}%
\Sigma=\operatorname*{Cov}\left(  \left(
\begin{array}
[c]{c}%
A\\
B\\
C
\end{array}
\right)  \right)  =\left(
\begin{array}
[c]{ccc}%
1 & \sigma & \rho\\
\sigma & 1 & \rho\\
\rho & \rho & 1
\end{array}
\right)  , &  & 1-\sigma^{2}>0, &  & 1-2\rho^{2}+\sigma>0.
\end{array}
\]
We wish to evaluate the probability that $AC<x$ and $BC<y$. On the one hand,
the joint density of $(AC,BC)$ is \cite{Mi1}
\[
f(x,y)=\frac{1}{2\pi}\frac{\exp\left(  \dfrac{1}{\xi}\left(  \rho x+\rho
y-\sqrt{\eta}\sqrt{(1-\rho^{2})x^{2}-2(-\rho^{2}+\sigma)xy+(1-\rho^{2})y^{2}%
}\right)  \right)  }{\sqrt{(1-\rho^{2})x^{2}-2(-\rho^{2}+\sigma)xy+(1-\rho
^{2})y^{2}}}%
\]
where $\xi=1-2\rho^{2}+\sigma$, $\eta=(1+\sigma)/(1-\sigma)$. On the other
hand, $(A^{2},2AB,2AC,B^{2},2BC,C^{2})$ is pseudo-Wishart distributed with $1$
degree of freedom and thus $(AC,BC)$ has characteristic function \cite{And}
\begin{align*}
F(u,v)  & =\left(  \det\left[  \left(
\begin{array}
[c]{ccc}%
1 & 0 & 0\\
0 & 1 & 0\\
0 & 0 & 1
\end{array}
\right)  -2i\left(
\begin{array}
[c]{ccc}%
1 & \sigma & \rho\\
\sigma & 1 & \rho\\
\rho & \rho & 1
\end{array}
\right)  \left(
\begin{array}
[c]{ccc}%
0 & 0 & u/2\\
0 & 0 & v/2\\
u/2 & v/2 & 0
\end{array}
\right)  \right]  \right)  ^{-1/2}\\
& =\frac{1}{\sqrt{1+(1-\rho^{2})u^{2}-2i\rho u+(1-\rho^{2})v^{2}-2i\rho
v+2(-\rho^{2}+\sigma)uv}}.
\end{align*}
Our first task is to confirm that $f(x,y)$ and $F(u,v)$ are indeed a Fourier
transform pair.

Now let $(A_{1},B_{1},C_{1})$, $(A_{2},B_{2},C_{2})$, $\ldots$, $(A_{n}%
,B_{n},C_{n})$ be a random sample from $N(0,\Sigma)$ and define sample
covariances
\[%
\begin{array}
[c]{ccc}%
\hat\gamma_{A,C}=%
{\displaystyle\sum\limits_{j=1}^{n}}
A_{j}C_{j}, &  & \hat\gamma_{B,C}=%
{\displaystyle\sum\limits_{j=1}^{n}}
B_{j}C_{j}.
\end{array}
\]
The marginal density for $\hat\gamma_{A,C}$ is well-known \cite{PJE, WB, Hir,
MBR, Mi2, Pr1, Spr, JO}:
\[
\frac{|x/2|^{_{(n-1)/2}}}{\sqrt{\pi}\sqrt{1-\rho^{2}}\Gamma\left(  n/2\right)
}\exp\left(  \frac{\rho x}{1-\rho^{2}}\right)  K_{(n-1)/2}\left(  \frac
{|x|}{1-\rho^{2}}\right)
\]
and likewise for $\hat\gamma_{B,C}$, where $K_{_{(n-1)/2}}(\theta)$ is the
modified Bessel function of the second kind. The joint density $f_{n}(x,y)$,
however, is not fully understood \cite{Pr2} even though
\[
F_{n}(u,v)=\frac1{\left[  1+(1-\rho^{2})u^{2}-2i\rho u+(1-\rho^{2}%
)v^{2}-2i\rho v+2(-\rho^{2}+\sigma)uv\right]  ^{n/2}}
\]
is comparatively simple. We shall determine $f_{n}(x,y)$ for $n=2,3,4$ and
then for arbitrary $n$, closing an evidently open issue. A\ special case
involving $f_{3}(x,y)$ was examined in \cite{F1}, for which $\rho=\sigma=1/2$,
to answer a geometric probability question. Our discussion generalizes this
earlier work.

\section{Fourier Transform Pair}

Our objective is to evaluate the integral:
\[
F(u,v)=\dfrac1{2\pi}%
{\displaystyle\int\limits_{-\infty}^{\infty}}
{\displaystyle\int\limits_{-\infty}^{\infty}}
\dfrac{\exp\left(  \dfrac1\xi\left(  (\rho+i\xi u)x+(\rho+i\xi v)y-\sqrt{\eta
}\sqrt{ax^{2}-2bxy+ay^{2}}\right)  \right)  }{\sqrt{ax^{2}-bxy+ay^{2}}}dx\,dy
\]
where $a=1-\rho^{2}$, $b=-\rho^{2}+\sigma$. Let
\[%
\begin{array}
[c]{ccc}%
x=\dfrac r{\sqrt{2}}\left(  \lambda\cos(\theta)-\kappa\sin(\theta)\right)  , &
& y=\dfrac r{\sqrt{2}}\left(  \lambda\cos(\theta)+\kappa\sin(\theta)\right)
\end{array}
\]
then $ax^{2}-bxy+by^{2}=r^{2}$ under the requirement that
\[%
\begin{array}
[c]{ccc}%
\lambda=\dfrac1{\sqrt{1-\sigma}}, &  & \kappa=\dfrac1{\sqrt{1-2\rho^{2}%
+\sigma}}.
\end{array}
\]
The Jacobian determinant is
\[
\left(  \dfrac1{\sqrt{2}}\right)  ^{2}\left|
\begin{array}
[c]{ccc}%
\lambda\cos(\theta)-\kappa\sin(\theta) &  & r\left(  -\lambda\sin
(\theta)-\kappa\cos(\theta)\right) \\
\lambda\cos(\theta)+\kappa\sin(\theta) &  & r\left(  -\lambda\sin
(\theta)+\kappa\cos(\theta)\right)
\end{array}
\right|  =\lambda\kappa r
\]
hence $dx\,dy=\lambda\kappa r\,dr\,d\theta$. We obtain
\begin{align*}
F(u,v)  & =\dfrac{\lambda\kappa}{2\pi}%
{\displaystyle\int\limits_{0}^{2\pi}}
{\displaystyle\int\limits_{0}^{\infty}}
\exp\left(  \frac1\xi\dfrac r{\sqrt{2}}\left[  \left(  \rho+i\xi u\right)
\left(  \lambda\cos(\theta)-\kappa\sin(\theta)\right)  \right.  \right. \\
& \;\;\;\;\;\;\;\;\;\;\;\;\;\left.  \left.  +\left(  \rho+i\xi v\right)
\left(  \lambda\cos(\theta)+\kappa\sin(\theta)\right)  -\sqrt{2\eta}\right]
\right)  dr\,d\theta\\
& =\dfrac{\lambda\kappa}{2\pi}%
{\displaystyle\int\limits_{0}^{2\pi}}
{\displaystyle\int\limits_{0}^{\infty}}
\exp\left(  \frac1\xi\dfrac r{\sqrt{2}}\left[  2\lambda\rho\cos(\theta
)+i\lambda\xi(u+v)\cos(\theta)-i\kappa\xi(u-v)\sin(\theta)-\sqrt{2\eta
}\right]  \right)  dr\,d\theta\\
\  & =-\dfrac{\lambda\kappa\xi}{\sqrt{2}\pi}%
{\displaystyle\int\limits_{0}^{2\pi}}
\frac1{2\lambda\rho\cos(\theta)+i\lambda\xi(u+v)\cos(\theta)-i\kappa
\xi(u-v)\sin(\theta)-\sqrt{2\eta}}d\theta
\end{align*}
since
\begin{align*}
\operatorname*{Re}\left[  2\lambda\rho\cos(\theta)+i\lambda\xi(u+v)\cos
(\theta)-i\kappa\xi(u-v)\sin(\theta)-\sqrt{2\eta}\right]   & =2\lambda\rho
\cos(\theta)-\sqrt{2\eta}\\
& \leq2\lambda\rho-\sqrt{2\eta}%
\end{align*}
and
\[
4\lambda^{2}\rho^{2}-2\eta=4\frac1{1-\sigma}\rho^{2}-2\frac{1+\sigma}%
{1-\sigma}=-\frac2{1-\sigma}\left(  1-2\rho^{2}+\sigma\right)  <0.
\]
Let $z=\exp(i\,\theta)$, then $d\theta=-i\,dz/z$ and
\begin{align*}
F(u,v)  & =\dfrac{i\lambda\kappa\xi}{\sqrt{2}\pi}%
{\displaystyle\int\nolimits_{C}}
\frac1{2\lambda\rho\frac12(z+\frac1z)+i\lambda\xi(u+v)\frac12(z+\frac
1z)-i\kappa\xi(u-v)\frac1{2i}(z-\frac1z)-\sqrt{2\eta}}\frac{dz}z\\
\  & =\dfrac{\sqrt{2}i\lambda\kappa\xi}\pi%
{\displaystyle\int\nolimits_{C}}
\frac1{2\lambda\rho(z^{2}+1)+i\lambda\xi(u+v)(z^{2}+1)-\kappa\xi
(u-v)(z^{2}-1)-2\sqrt{2\eta}z}dz.
\end{align*}
The denominator of the integrand can be rewritten as
\[
(2\lambda\rho-\kappa\xi u+\kappa\xi v+i\lambda\xi u+i\lambda\xi v)z^{2}%
-2\sqrt{2\eta}z+(2\lambda\rho+\kappa\xi u-\kappa\xi v+i\lambda\xi
u+i\lambda\xi v).
\]
Let $\delta=2\rho/\xi$. The two poles $z_{\text{pos}}$, $z_{\text{neg}}$ of
the integrand are
\begin{align*}
& \frac{2\sqrt{2\eta}\pm\sqrt{8\eta-4(2\lambda\rho-\kappa\xi u+\kappa\xi
v+i\lambda\xi u+i\lambda\xi v)(2\lambda\rho+\kappa\xi u-\kappa\xi
v+i\lambda\xi u+i\lambda\xi v)}}{2(2\lambda\rho-\kappa\xi u+\kappa\xi
v+i\lambda\xi u+i\lambda\xi v)}\\
& =\frac{\sqrt{2\eta}\pm\sqrt{\lambda^{2}\xi^{2}(u+v-i\delta)^{2}+\kappa
^{2}\xi^{2}(-u+v)^{2}+2\eta}}{2\lambda\rho-\kappa\xi u+\kappa\xi v+i\lambda\xi
u+i\lambda\xi v}%
\end{align*}
and $z_{\text{neg}}$ is always inside $C$, $z_{\text{pos}}$ is always outside.
Clearly $z_{\text{neg}}$ is a pole of order $1$ and the associated residue is
\begin{align*}
& \lim_{z\rightarrow z_{\text{neg}}}\frac1{(2\lambda\rho-\kappa\xi u+\kappa\xi
v+i\lambda\xi u+i\lambda\xi v)(z-z_{\text{pos}})}\\
& =-\frac12\frac1{\sqrt{\lambda^{2}\xi^{2}(u+v-i\delta)^{2}+\kappa^{2}\xi
^{2}(-u+v)^{2}+2\eta}};
\end{align*}
multiplying by $(2\pi i)(\sqrt{2}i\lambda\kappa\xi/\pi)$ gives
\begin{align*}
& \frac{\sqrt{2}\lambda\kappa\xi}{\sqrt{\lambda^{2}\xi^{2}(u+v-i\delta
)^{2}+\kappa^{2}\xi^{2}(-u+v)^{2}+2\eta}}\\
& =\frac{\sqrt{2}\lambda\kappa}{\sqrt{\lambda^{2}(u+v-i\delta)^{2}+\kappa
^{2}(-u+v)^{2}+\alpha^{2}}}%
\end{align*}
where $\alpha=\sqrt{2\eta}/\xi$. This is the most useful expression for our
purposes. The original formula for $F(u,v)$ can now be confirmed.

For example, if $\rho=\sigma=1/2$, we have $\xi=1$, $\eta=3$, $\lambda
=\sqrt{2}$, $\kappa=1$, $\delta=1$, $\alpha=\sqrt{6}$. In this special case,
our expression for $F(u,v)$ becomes
\[
\frac2{\sqrt{2(u+v-i)^{2}+(-u+v)^{2}+6}}=\frac1{\sqrt{(\frac u2-i)(3\frac
u2+i)+(\frac v2-i)(3\frac v2+i)+2\frac u2\frac v2-1}}
\]
and this is consistent with \cite{F1}.

We can also argue in reverse, that is, evaluating instead the integral:
\[
f(x,y)=\frac{\sqrt{2}\lambda\kappa}{(2\pi)^{2}}%
{\displaystyle\int\limits_{-\infty}^{\infty}}
{\displaystyle\int\limits_{-\infty}^{\infty}}
\frac{\exp(-iux-ivy)}{\sqrt{\lambda^{2}(u+v-i\delta)^{2}+\kappa^{2}%
(-u+v)^{2}+\alpha^{2}}}dv\,du.
\]
Let $s=u+v$, $t=-u+v$, then $u=(s-t)/2$, $v=(s+t)/2$ and
\begin{align*}
-iux-ivy  & =-i\frac s2x+i\frac t2x-i\frac s2y-i\frac t2y\\
\  & =-i\frac{x+y}2s-i\frac{-x+y}2t;
\end{align*}
hence
\[
f(x,y)=\frac{\lambda\kappa}{\sqrt{2}(2\pi)^{2}}%
{\displaystyle\int\limits_{-\infty}^{\infty}}
{\displaystyle\int\limits_{-\infty}^{\infty}}
\frac{\exp\left(  -i\frac{x+y}2s-i\frac{-x+y}2t\right)  }{\sqrt{\lambda
^{2}(s-i\delta)^{2}+\kappa^{2}t^{2}+\alpha^{2}}}dt\,ds
\]
since the Jacobian determinant is $1/2$. The inner integral is \cite{GR}
\[%
{\displaystyle\int\limits_{-\infty}^{\infty}}
\frac{\exp\left(  -i\frac{-x+y}2t\right)  }{\sqrt{\lambda^{2}(s-i\delta
)^{2}+\kappa^{2}t^{2}+\alpha^{2}}}dt=\frac2\kappa K_{0}\left(  \frac
1\kappa\left|  \frac{-x+y}2\right|  \sqrt{\lambda^{2}(s-i\delta)^{2}%
+\alpha^{2}}\right)  .
\]
Define
\[
\gamma=\frac1\kappa\left|  \frac{-x+y}2\right|  ,
\]
then the outer integral is \cite{GR}
\begin{align*}
&
{\displaystyle\int\limits_{-\infty}^{\infty}}
K_{0}\left(  \gamma\sqrt{\lambda^{2}(s-i\delta)^{2}+\alpha^{2}}\right)
\exp\left(  -i\frac{x+y}2s\right)  ds\\
\  & =\pi\exp\left(  \frac{x+y}2\delta\right)  \frac{\exp\left(  -\frac
\alpha\lambda\sqrt{\gamma^{2}\lambda^{2}+\left(  \frac{x+y}2\right)  ^{2}%
}\right)  }{\sqrt{\gamma^{2}\lambda^{2}+\left(  \frac{x+y}2\right)  ^{2}}}\\
\  & =2\pi\kappa\exp\left(  \frac{x+y}2\delta\right)  \frac{\exp\left(
-\frac\alpha{2\lambda\kappa}\sqrt{\kappa^{2}\left(  x+y\right)  ^{2}%
+\lambda^{2}\left(  x-y\right)  ^{2}}\right)  }{\sqrt{\kappa^{2}\left(
x+y\right)  ^{2}+\lambda^{2}\left(  x-y\right)  ^{2}}}.
\end{align*}
Multiplying by $(\lambda\kappa/(\sqrt{2}(2\pi)^{2}))\cdot(2/\kappa)$, we
obtain
\[
\frac{\lambda\kappa}{\sqrt{2}\pi}\exp\left(  \frac{x+y}2\delta\right)
\frac{\exp\left(  -\frac\alpha{2\lambda\kappa}\sqrt{\kappa^{2}\left(
x+y\right)  ^{2}+\lambda^{2}\left(  x-y\right)  ^{2}}\right)  }{\sqrt
{\kappa^{2}\left(  x+y\right)  ^{2}+\lambda^{2}\left(  x-y\right)  ^{2}}}
\]
and the original formula for $f(x,y)$ can now be confirmed. For example, if
$\rho=\sigma=1/2$, the probability that both $AC>0$ and $BC>0$ is $1/2=0.5$.

\section{Sample Size $n=2$}

We wish to evaluate
\[
f_{2}(x,y)=\frac{\lambda^{2}\kappa^{2}}{(2\pi)^{2}}%
{\displaystyle\int\limits_{-\infty}^{\infty}}
{\displaystyle\int\limits_{-\infty}^{\infty}}
\frac{\exp\left(  -i\frac{x+y}2s-i\frac{-x+y}2t\right)  }{\lambda
^{2}(s-i\delta)^{2}+\kappa^{2}t^{2}+\alpha^{2}}dt\,ds.
\]
The inner integral is \cite{GR}
\[%
{\displaystyle\int\limits_{-\infty}^{\infty}}
\frac{\exp\left(  -i\frac{-x+y}2t\right)  }{\lambda^{2}(s-i\delta)^{2}%
+\kappa^{2}t^{2}+\alpha^{2}}dt=\frac\pi\kappa\frac{\exp\left(  -\frac
1\kappa\left|  \frac{-x+y}2\right|  \sqrt{\lambda^{2}(s-i\delta)^{2}%
+\alpha^{2}}\right)  }{\sqrt{\lambda^{2}(s-i\delta)^{2}+\alpha^{2}}}
\]
and the outer integral is \cite{GR}
\begin{align*}
&
{\displaystyle\int\limits_{-\infty}^{\infty}}
\frac{\exp\left(  -\gamma\sqrt{\lambda^{2}(s-i\delta)^{2}+\alpha^{2}}\right)
}{\sqrt{\lambda^{2}(s-i\delta)^{2}+\alpha^{2}}}\exp\left(  -i\frac
{x+y}2s\right)  ds\\
\  & =\frac2\lambda\exp\left(  \frac{x+y}2\delta\right)  K_{0}\left(
\frac\alpha\lambda\sqrt{\gamma^{2}\lambda^{2}+\left(  \frac{x+y}2\right)
^{2}}\right) \\
\  & =\frac2\lambda\exp\left(  \frac{x+y}2\delta\right)  K_{0}\left(
\frac\alpha{2\lambda\kappa}\sqrt{\kappa^{2}\left(  x+y\right)  ^{2}%
+\lambda^{2}\left(  x-y\right)  ^{2}}\right)  .
\end{align*}
Multiplying by $(\lambda^{2}\kappa^{2}/(2\pi)^{2}))\cdot(\pi/\kappa)$, we
obtain that $f_{2}(x,y)$ is
\[
\ \frac{\lambda\kappa}{2\pi}\exp\left(  \frac{x+y}2\delta\right)  K_{0}\left(
\frac\alpha{2\lambda\kappa}\sqrt{\kappa^{2}\left(  x+y\right)  ^{2}%
+\lambda^{2}\left(  x-y\right)  ^{2}}\right)  .
\]
For example, if $\rho=\sigma=1/2$, the probability that both $\hat\gamma
_{A,C}>0$ and $\hat\gamma_{B,C}>0$ is $0.608173447....$

\section{Sample Size $n=3$}

We wish to evaluate
\[
f_{3}(x,y)=\frac{\sqrt{2}\lambda^{3}\kappa^{3}}{(2\pi)^{2}}%
{\displaystyle\int\limits_{-\infty}^{\infty}}
{\displaystyle\int\limits_{-\infty}^{\infty}}
\frac{\exp\left(  -i\frac{x+y}2s-i\frac{-x+y}2t\right)  }{\left(  \lambda
^{2}(s-i\delta)^{2}+\kappa^{2}t^{2}+\alpha^{2}\right)  ^{3/2}}dt\,ds.
\]
The inner integral is
\begin{align*}%
{\displaystyle\int\limits_{-\infty}^{\infty}}
\frac{\exp\left(  -i\frac{-x+y}2t\right)  }{\left(  \lambda^{2}(s-i\delta
)^{2}+\kappa^{2}t^{2}+\alpha^{2}\right)  ^{3/2}}dt  & =-\frac1{\lambda
^{2}(s-i\delta)}\frac d{ds}%
{\displaystyle\int\limits_{-\infty}^{\infty}}
\frac{\exp\left(  -i\frac{-x+y}2t\right)  }{\sqrt{\lambda^{2}(s-i\delta
)^{2}+\kappa^{2}t^{2}+\alpha^{2}}}dt\\
& =-\frac1{\lambda^{2}(s-i\delta)}\frac d{ds}\frac2\kappa K_{0}\left(
\frac1\kappa\left|  \frac{-x+y}2\right|  \sqrt{\lambda^{2}(s-i\delta
)^{2}+\alpha^{2}}\right) \\
& =\frac2{\kappa^{2}}\left|  \frac{-x+y}2\right|  \frac{K_{1}\left(
\frac1\kappa\left|  \frac{-x+y}2\right|  \sqrt{\lambda^{2}(s-i\delta
)^{2}+\alpha^{2}}\right)  }{\sqrt{\lambda^{2}(s-i\delta)^{2}+\alpha^{2}}}%
\end{align*}
and the outer integral is \cite{GR}
\begin{align*}
&
{\displaystyle\int\limits_{-\infty}^{\infty}}
\frac{K_{1}\left(  \gamma\sqrt{\lambda^{2}(s-i\delta)^{2}+\alpha^{2}}\right)
}{\sqrt{\lambda^{2}(s-i\delta)^{2}+\alpha^{2}}}\exp\left(  -i\frac
{x+y}2s\right)  ds\\
\  & =\frac\pi{\alpha\gamma\lambda}\exp\left(  \frac{x+y}2\delta\right)
\exp\left(  -\frac\alpha\lambda\sqrt{\gamma^{2}\lambda^{2}+\left(  \frac
{x+y}2\right)  ^{2}}\right) \\
\  & =\frac\pi{\alpha\gamma\lambda}\exp\left(  \frac{x+y}2\delta\right)
\exp\left(  -\frac\alpha{2\lambda\kappa}\sqrt{\kappa^{2}\left(  x+y\right)
^{2}+\lambda^{2}\left(  x-y\right)  ^{2}}\right)  .
\end{align*}
Multiplying by $(\sqrt{2}\lambda^{3}\kappa^{3}/(2\pi)^{2})\cdot(2\gamma
/\kappa)$, we obtain that $f_{3}(x,y)$ is
\begin{align*}
& \ \ \ \frac{\lambda^{2}\kappa^{2}}{\sqrt{2}\pi\alpha}\exp\left(  \frac
{x+y}2\delta\right)  \exp\left(  -\frac\alpha{2\lambda\kappa}\sqrt{\kappa
^{2}\left(  x+y\right)  ^{2}+\lambda^{2}\left(  x-y\right)  ^{2}}\right) \\
\  & =\frac1{2\pi\sqrt{1-\sigma^{2}}}\exp\left(  \dfrac1\xi\left(  \rho x+\rho
y-\sqrt{\eta}\sqrt{(1-\rho^{2})x^{2}-2(-\rho^{2}+\sigma)xy+(1-\rho^{2})y^{2}%
}\right)  \right)
\end{align*}
which is remarkably simple. For example, if $\rho=\sigma=1/2$, the probability
that both $\hat\gamma_{A,C}>0$ and $\hat\gamma_{B,C}>0$ is $0.6837762984...$.

\section{Sample Size $n=4$}

We wish to evaluate
\[
f_{4}(x,y)=\frac{2\lambda^{4}\kappa^{4}}{(2\pi)^{2}}%
{\displaystyle\int\limits_{-\infty}^{\infty}}
{\displaystyle\int\limits_{-\infty}^{\infty}}
\frac{\exp\left(  -i\frac{x+y}2s-i\frac{-x+y}2t\right)  }{\left(  \lambda
^{2}(s-i\delta)^{2}+\kappa^{2}t^{2}+\alpha^{2}\right)  ^{2}}dt\,ds.
\]
The inner integral is \cite{GR}
\begin{align*}%
{\displaystyle\int\limits_{-\infty}^{\infty}}
\frac{\exp\left(  -i\frac{-x+y}2t\right)  }{\left(  \lambda^{2}(s-i\delta
)^{2}+\kappa^{2}t^{2}+\alpha^{2}\right)  ^{2}}dt  & =\frac\pi{2\kappa^{2}%
}\left|  \frac{-x+y}2\right|  \exp\left(  -\frac1\kappa\left|  \frac
{-x+y}2\right|  \sqrt{\lambda^{2}(s-i\delta)^{2}+\alpha^{2}}\right)  \cdot\\
& \ \ \ \ \left\{  \frac1{\frac1\kappa\left|  \frac{-x+y}2\right|  \left(
\lambda^{2}(s-i\delta)^{2}+\alpha^{2}\right)  ^{3/2}}+\frac1{\lambda
^{2}(s-i\delta)^{2}+\alpha^{2}}\right\}
\end{align*}
and the outer integral is \cite{GR}
\begin{align*}
&
{\displaystyle\int\limits_{-\infty}^{\infty}}
\left\{  \frac1{\gamma\left(  \lambda^{2}(s-i\delta)^{2}+\alpha^{2}\right)
^{3/2}}+\frac1{\lambda^{2}(s-i\delta)^{2}+\alpha^{2}}\right\}  \exp\left(
-\gamma\sqrt{\lambda^{2}(s-i\delta)^{2}+\alpha^{2}}\right)  \exp\left(
-i\frac{x+y}2s\right)  ds\\
\  & =\frac2{\alpha\gamma\lambda^{2}}\exp\left(  \frac{x+y}2\delta\right)
\sqrt{\gamma^{2}\lambda^{2}+\left(  \frac{x+y}2\right)  ^{2}}K_{1}\left(
\frac\alpha\lambda\sqrt{\gamma^{2}\lambda^{2}+\left(  \frac{x+y}2\right)
^{2}}\right) \\
\  & =\frac1{\alpha\gamma\lambda^{2}\kappa}\exp\left(  \frac{x+y}%
2\delta\right)  \sqrt{\kappa^{2}\left(  x+y\right)  ^{2}+\lambda^{2}\left(
x-y\right)  ^{2}}K_{1}\left(  \frac\alpha{2\lambda\kappa}\sqrt{\kappa
^{2}\left(  x+y\right)  ^{2}+\lambda^{2}\left(  x-y\right)  ^{2}}\right)  .
\end{align*}
Multiplying by $(2\lambda^{4}\kappa^{4}/(2\pi)^{2})\cdot(\pi\gamma/(2\kappa
))$, we obtain that $f_{4}(x,y)$ is
\[
\ \ \frac{\lambda^{2}\kappa^{2}}{4\pi\alpha}\exp\left(  \frac{x+y}%
2\delta\right)  \sqrt{\kappa^{2}\left(  x+y\right)  ^{2}+\lambda^{2}\left(
x-y\right)  ^{2}}K_{1}\left(  \frac\alpha{2\lambda\kappa}\sqrt{\kappa
^{2}\left(  x+y\right)  ^{2}+\lambda^{2}\left(  x-y\right)  ^{2}}\right)  .
\]
For example, if $\rho=\sigma=1/2$, the probability that both $\hat\gamma
_{A,C}>0$ and $\hat\gamma_{B,C}>0$ is $0.7409625593....$

\section{Formula for Arbitrary $n$}

We wish to evaluate
\[
f_{n}(x,y)=\frac{2^{(n-2)/2}\lambda^{n}\kappa^{n}}{(2\pi)^{2}}%
{\displaystyle\int\limits_{-\infty}^{\infty}}
{\displaystyle\int\limits_{-\infty}^{\infty}}
\frac{\exp\left(  -i\frac{x+y}2s-i\frac{-x+y}2t\right)  }{\left(  \lambda
^{2}(s-i\delta)^{2}+\kappa^{2}t^{2}+\alpha^{2}\right)  ^{n/2}}dt\,ds.
\]
The inner integral is \cite{Mc, Wa, AS}
\begin{align*}%
{\displaystyle\int\limits_{-\infty}^{\infty}}
\frac{\exp\left(  -i\frac{-x+y}2t\right)  }{\left(  \lambda^{2}(s-i\delta
)^{2}+\kappa^{2}t^{2}+\alpha^{2}\right)  ^{n/2}}dt  & =\frac{\sqrt{2\pi}%
}{\Gamma(n/2)2^{(n-2)/2}}\frac1{\kappa^{(n+1)/2}}\left|  \frac{-x+y}2\right|
^{(n-1)/2}\frac1{\left(  \lambda^{2}(s-i\delta)^{2}+\alpha^{2}\right)
^{(n-1)/4}}\\
& \ \cdot K_{(n-1)/2}\left(  \frac1\kappa\left|  \frac{-x+y}2\right|
\sqrt{\lambda^{2}(s-i\delta)^{2}+\alpha^{2}}\right)
\end{align*}
and the outer integral is \cite{GR}
\begin{align*}
& \ \ \
{\displaystyle\int\limits_{-\infty}^{\infty}}
\frac{K_{(n-1)/2}\left(  \gamma\sqrt{\lambda^{2}(s-i\delta)^{2}+\alpha^{2}%
}\right)  }{\left(  \lambda^{2}(s-i\delta)^{2}+\alpha^{2}\right)  ^{(n-1)/4}%
}\exp\left(  -i\frac{x+y}2s\right)  ds\\
\  & =\sqrt{\frac\pi2}\frac2{\alpha^{(n-2)/2}\gamma^{(n-1)/2}\lambda^{n/2}%
}\exp\left(  \frac{x+y}2\delta\right)  \left(  \gamma^{2}\lambda^{2}+\left(
\frac{x+y}2\right)  ^{2}\right)  ^{(n-2)/4}\\
& \ \cdot K_{(n-2)/2}\left(  \frac\alpha\lambda\sqrt{\gamma^{2}\lambda
^{2}+\left(  \frac{x+y}2\right)  ^{2}}\right) \\
\  & =\sqrt{\frac\pi2}\frac1{2^{(n-4)/2}\alpha^{(n-2)/2}\gamma^{(n-1)/2}%
\lambda^{n/2}\kappa^{(n-2)/2}}\exp\left(  \frac{x+y}2\delta\right)  \left(
\kappa^{2}\left(  x+y\right)  ^{2}+\lambda^{2}\left(  x-y\right)  ^{2}\right)
^{(n-2)/4}\\
& \ \ \ \cdot K_{(n-2)/2}\left(  \frac\alpha{2\lambda\kappa}\sqrt{\kappa
^{2}\left(  x+y\right)  ^{2}+\lambda^{2}\left(  x-y\right)  ^{2}}\right)  .
\end{align*}
Multiplying by
\[
\frac{2^{(n-2)/2}\lambda^{n}\kappa^{n}}{(2\pi)^{2}}\cdot\frac{\sqrt{2\pi}%
}{\Gamma(n/2)2^{(n-2)/2}}\frac1\kappa\gamma^{(n-1)/2},
\]
we obtain that $f_{n}(x,y)$ is
\begin{align*}
& \ \ \ \ \ \frac1{\Gamma(n/2)}\frac{\lambda^{n/2}\kappa^{n/2}}{2^{n/2}%
\pi\alpha^{(n-2)/2}}\exp\left(  \frac{x+y}2\delta\right)  \left(  \kappa
^{2}\left(  x+y\right)  ^{2}+\lambda^{2}\left(  x-y\right)  ^{2}\right)
^{(n-2)/4}\\
& \ \ \ \cdot K_{(n-2)/2}\left(  \frac\alpha{2\lambda\kappa}\sqrt{\kappa
^{2}\left(  x+y\right)  ^{2}+\lambda^{2}\left(  x-y\right)  ^{2}}\right)
\end{align*}
which indeed generalizes the cases $n=1,2,3,4$ worked earlier. This result is
asymptotically consistent with the following outcome of the Central Limit
Theorem:
\[
\left(
\begin{array}
[c]{c}%
(\hat\gamma_{A,C}-n\rho)/\sqrt{n}\\
(\hat\gamma_{B,C}-n\rho)/\sqrt{n}%
\end{array}
\right)  \sim N\left(  \left(
\begin{array}
[c]{c}%
0\\
0
\end{array}
\right)  ,\left(
\begin{array}
[c]{cc}%
\rho^{2}+1 & \rho^{2}+\sigma\\
\rho^{2}+\sigma & \rho^{2}+1
\end{array}
\right)  \right)
\]
as $n\rightarrow\infty$.

\section{Closing Words}

Testing the hypothesis $H_{0}:\gamma_{A,C}=\gamma_{B,C}$ can be done by
examining the difference
\[
\hat\gamma_{A,C}-\hat\gamma_{B,C}=\hat\gamma_{A-B,C}=%
{\displaystyle\sum\limits_{j=1}^{n}}
(A_{j}-B_{j})C_{j}.
\]
If $H_{0}$ is true, then
\[
\left(
\begin{array}
[c]{c}%
A-B\\
C
\end{array}
\right)  \sim N\left(  \left(
\begin{array}
[c]{c}%
0\\
0
\end{array}
\right)  ,\left(
\begin{array}
[c]{cc}%
2-2\sigma & 0\\
0 & 1
\end{array}
\right)  \right)
\]
and is independent of $\rho$; further, the density of $\hat\gamma_{A-B,C}$ is
\[
\frac{|x/2|^{_{(n-1)/2}}}{\sqrt{\pi}\left(  2-2\sigma\right)  ^{(n+1)/4}%
\Gamma\left(  n/2\right)  }K_{(n-1)/2}\left(  \frac{|x|}{\sqrt{2-2\sigma}%
}\right)
\]
via known results \cite{Mi2} on Gaussian inner products. A\ considerable
literature exists on the harder problem of testing $\tilde H_{0}:\rho
_{A,C}=\rho_{B,C}$ where variances are unknown and underlying distributions
might not be normal \cite{X1, X2, X3, X4, X5, X6, X7, X8, X9, X10, X11, X12,
X13, X14, X15, X16, X17, X18, X19, X20, X21, X22, X23, X24, X25, X26, X27,
X28, X29, X30, X31, X32, X33, X34}.

We have not attempted to evaluate the inverse Fourier transform of%

\[
G(u,v)=\left(  \det\left[  \left(
\begin{array}
[c]{ccc}%
1 & 0 & 0\\
0 & 1 & 0\\
0 & 0 & 1
\end{array}
\right)  -2i\left(
\begin{array}
[c]{ccc}%
1 & \sigma & \rho\\
\sigma & 1 & \rho\\
\rho & \rho & 1
\end{array}
\right)  \left(
\begin{array}
[c]{ccc}%
0 & u/2 & v/2\\
u/2 & 0 & w/2\\
v/2 & w/2 & 0
\end{array}
\right)  \right]  \right)  ^{-1/2}.
\]
There is no analog of Miller's result \cite{Mi1}, as far as we know, giving a
joint density $g(x,y,z)$ for $(AB,AC,BC)$. Hence no distributional insight on
$(\hat\gamma_{A,B},\hat\gamma_{A,C},\hat\gamma_{B,C})$ seems to be available
here. Interestingly, a formula for a\ joint density for $(\hat\rho_{A,B}%
,\hat\rho_{A,C},\hat\rho_{B,C})$ is outlined in \cite{Fs, Be} -- evidently a
sample size $>4$ is presumed -- and details still need to come together.

\section{Acknowledgement}

I am grateful to Robert Israel for a helpful discussion about residue
calculus. Much more relevant material can be found at \cite{F2, F3}, including
experimental computer runs that aided theoretical discussion here.

\end{document}